
\input amstex
\documentstyle{amsppt}
\rightheadtext{On $CR$ mappings and separate algebraicity.}
\nopagenumbers
\topmatter
\title
On $CR$ mappings between algebraic Cauchy--Riemann manifolds
and separate algebraicity for holomorphic functions.
\endtitle
\author
Ruslan Sharipov and Alexander Sukhov
\endauthor
\abstract
    We prove the algebraicity of smooth $CR$-mappings between
algebraic Cauchy--Riemann manifolds. A generalization of separate
algebraicity principle is established.
\endabstract
\address
Department of Mathematics, Bashkir State University,
Frunze street 32, 450074 Ufa, Russia
\endaddress
\email
yavdat\@bgua.bashkiria.su
\endemail
\thanks
Paper is written under the financial support of International
Scientific Foundation of Soros, project \#RK4000, and Russian
Fund for Fundamental Researches, project \#93-01-00273.
\endthanks
\endtopmatter
\document
\head
1. Introduction
\endhead
The well-known Webster theorem \cite{W1} asserts that a local $CR$
diffeomorphism of a real algebraic Levi-nondegenerate hypersurface
in $\Bbb C^{\,n}$, for $n>1$, extends to an algebraic mapping of all
$\Bbb C{\,^n}$. In this paper we generalize this result for real
algebraic Cauchy--Riemann manifolds of higher codimensions (Theorem 1).
Our technique is based on the Webster reflection principle modified
in the spirit of \cite{Su} and a generalization of the classical
separate algebraicity theorem, where we replace the parallel lines
by families of algebraic curves (Theorem 2). We hope that this second
result is of self-interest.\par
     The paper is organized as follows. In section 2 we give the
precise definitions and statements of our results. In sections 3
and 4 we prove Theorem 1 provided Theorem 2 holds. Section 5 is
devoted to the proof of Theorem 2.\par
\head
2. The results.
\endhead
     Let $\Omega$ be a domain in $\Bbb C^{\,n}$. A closed subset
$M$ of $\Omega$ is called a generic real algebraic manifold of
codimension $d\geq 1$ if
$$
M=\{z\in\Omega: \rho_j(z,\bar z)=0, j=1,\dots,d\}\tag2.1
$$
where $\rho_j$ are real polynomials and $\bar\partial\rho_1\land
\ldots\land\bar\partial\rho_d\neq 0$ in $\Omega$. They are called
the defining functions of $M$. We denote by $T_pM$ and $T_p^cM$
the real and complex tangent spaces of $M$ at the point $p\in M$
(recall that $T^c_pM=T_pM\cap J(T_pM)$, where $J$ is the complex
structure operator in $\Bbb C^{\,n}$). This is well known that for
a generic manifold the complex dimension of $T_p^cM$ does not
depend on $p\in M$ and is equal to $n-d$; it is called the
$CR$-dimension of $M$.\par
     We denote by $H_p(\rho_j,u,v)$ the value of the Levi form
(complex hessian) of the function $\rho_j$ on vectors $u$, $v$
at the point $p\in M$, i.e.
$$
H_p(\rho_j,u,v)=\sum^n_{\nu,\mu=1}\frac{\partial^2\rho_j}
{\partial z_\nu\,\partial\bar z_\mu}(p)\,u_\nu\bar v_\mu
\tag2.2
$$\par
    The Levi cone (at $p\in M$) of the manifold $M$ of the
form \thetag{2.1} is said to be the convex hull of the set
$$
\{\alpha=(\alpha_1,\dots,\alpha_d)\in\Bbb R^d: \alpha_j=
H_p(\rho_j,u,u), u\in T^c_pM\}
\tag2.3
$$\par
     If the Levi cone of $M$ has a non-empty interior in $\Bbb R^d$,
then we say that $M$ possess a non-degenerate Levi cone at $p$.
Evidently this condition does not depend on the choice of defining
functions and is invariant with respect to changes of coordinates.
\par
The vector valued hermitian form
$$
L_p(u,v)=(H_p(\rho_1,u,v),\dots,H_p(\rho_d,u,v))
$$
is called the Levi form of $M$ at $p$.\par
     Along with $M$ we consider a domain $\Omega'$ in $\Bbb C^{\,n'}$
and a generic real algebraic manifold $M'$ of codimension $d\,'\geq 1$
of the form
$$
M'=\{z'\in\Omega': \rho'_j(z',\bar z')=0, j=1,\dots,d\,'\}\tag2.4
$$
where $\rho'_j$ are real polynomials and $\bar\partial\rho'_1\land
\ldots\land\bar\partial\rho'_{d\,'}\neq 0$ in $\Omega'$. Fix a
hermitian scalar product in $\Bbb C^{\,n'}$ that is denoted by
$bigl<\cdot,\cdot\bigr>'$. Fix a point $p'\in M'$. With every
defining function $\rho'_j$ we associate the Levi operator
$L_{p'}^j: T^c_{p'}M'\longrightarrow T^c_{p'}M'$ determined by
$H_{p'}(\rho'_j,u,v)=\bigl<L^j_{p'}(u),\ v\bigr>'$ for all
$u,v\in T^c_{p'}M'$. Of course, this definition of the Levi operator
depends on the choice of hermitian scalar product. But as we shall
see further this dependence is unessential In what follows we
use the similar notation for the Levi operators of $M$ (without
primes).\par
     Let $U$ be an open connected subset of the manifold $M$.
The mapping $F: U\longrightarrow M'$ of the smoothness class $C^1$
is called the $CR$-mapping if for any point $p\in U$ the tangent map
$dF_p$ is $\Bbb C$\,-linear after restriction to the complex
tangent space $T^c_pM$ (in this case $dF_p(T^c_pM)\subset
T^c_{p'}M'$). We say that $F$ extends to an algebraic
mapping of all $\Bbb C^{\,n}$ if the graph of $F$ is a part of
$n$-dimensional complex algebraic manifold in $\Bbb C^{\,n+n'}$.
\par
    Now we can formulate the first main result of our paper. It is
given by the following theorem.
\proclaim{Theorem 1} Let $\Omega\in\Bbb C^{\,n}$ and $\Omega'\in
\Bbb C^{\,n'}$ be domains, $M\subset\Omega$ and $M'\subset
\Omega'$  be generic real algebraic manifolds of the form
\thetag{2.1} and \thetag{2.4} respectively, $M$ having the
non-degenerate Levi cone at some point $p\in M$. Suppose
$U\subset M$ is an open connected subset of $M$ containing $p$
and let $F: U\longrightarrow M'$ be a $CR$-mapping of class
$C^1$ satisfying the following condition
$$
\sum_{j=1}^{d\,'} L^j_{p'}(dF_p(T^c_pM))=T^c_{p'}M'
\text{, where }p'=F(p)
\tag2.5
$$
Then $F$ extends to an algebraic mapping of
the whole space $\Bbb C^{\,n}$.
\endproclaim
\par
     For the first note that \thetag{2.5} does not depend on the
choice of the hermitian scalar product in $\Bbb C^{\,n'}$ by
means of which we defined the Levi operators $L^j_{p'}$. If
$\tilde L^j_{p'}$ is defined by another scalar product, then
it is connected with $L^j_{p'}$ by the equality $\tilde L^j_{p'}=
A\,L^j_{p'}$, where $A$ is a non-degenerate $\Bbb C$\,-linear
operator in $T^c_{p'}M'$. Therefore \thetag{2.5} holds for
operators $\tilde L^j_{p'}$ as well.\par
Recall that the Levi form of $M$ is called non-degenerate if
$L_p(u,v)=0$ for any $v\in T^c_pM$ implies $u=0$ \cite{W2}.
\proclaim {Corollary}
Let $F:M\longrightarrow M'$ be a $CR$ diffeomorphism of class
$C^1$ between two real algebraic manifolds  in $\Bbb C^{\,n}$
with  non-degenerate Levi forms and non-degenerate Levi cones.
Then $F$ extends to an algebraic mapping on all $\Bbb C^{\,n}$.
\endproclaim
This assertion follows by theorem 1 quite similar to \cite{Su}
We emphasize that theorem 1 treats a considerably more general
situation, since $M$ and $M'$ are allowed to have different
$CR$ dimensions.
     Our proof of theorem 1 is based on the modification of the
Webster reflection principle \cite{W1}. The crucial technical
tool here is a "curved" version of the following classical separate
algebraicity principle.\par
\proclaim {Claim} Let $f(\bold z)=f(z_1,\dots,z_n)$ be a
function in some domain $D\subset C^{\,n}$. If
$f(\bold z)$ is algebraic in each separate variable $z_i$ for
any fixed values of other variables, then $f(\bold z)$ is an
algebraic function in $D$.
\endproclaim
     Proof of this classical theorem can be found in \cite{BM}.
Function $f(\bold z)$ in this assertion is algebraic
along the straight lines, parallel to the coordinate axis. The domain
$D$ foliates into $n$ families of such lines. In order to
generalize this result for our purposes we introduce
$n$ families of algebraic curves in $D$
$$
\aligned
&z_1=R^{(m)}_1(t_m,c^{(m)}_1,\dots,c^{(m)}_{n-1})\\
&\text{\hbox to 5cm{\leaders\hbox to 1em{\hss.\hss}\hfill}}\\
&z_n=R^{(m)}_n(t_m,c^{(m)}_1,\dots,c^{(m)}_{n-1})
\endaligned\tag2.6
$$
Here $m=1,\dots,n$ is the number of the family, $t_m$ is a parameter
on each particular curve of $m$-th family ($R^{(m)}_i$ depends on
$t_m$ algebraically). Parameters $c^{(m)}_1,\dots,c^{(m)}_{n-1}$
identify curves of $m$-th family. \par
\definition{Definition 1} The family of algebraic curves \thetag{2.6}
is called {\it algebraically depending on parameters} if each of the
defining functions $R^{(m)}_i$, $i=1,\dots,n$ in \thetag{2.6} is
algebraic in $c^{(m)}_1,\dots,c^{(m)}_{n-1}$.
\enddefinition
    In this case because of above classical separate algebraicity principle
functions $R^{(m)}_i$ in \thetag{2.6} are algebraic in whole
set of their arguments $t_m,c^{(m)}_1,\dots,c^{(m)}_{n-1}$.\par
\definition{Definition 2} The family of curves \thetag{2.6} is called
nonsingular in $D$ if curves of this family fills the whole domain $D$
and the mapping $R^{(m)}:\ (t_m,c^{(m)}_1,\dots,c^{(m)}_{n-1})\longrightarrow
\bold z$ is the local diffeomorphism.\par
\enddefinition
\definition{Definition 3} Let \thetag{2.6} define $n$ nonsingular
families of algebraic curves in $D$. Then we have $n$ tangent vectors
to the curves at each point of $D$. We shall say that $n$ families of
curves \thetag{2.6} are in general position if these vectors are
linearly independent at any point $\bold z\in D$.
\enddefinition
     Now we can state the "curved" separate algebraicity principle
which partially generalizes the classical one and form the second
main result of our paper.
\proclaim{Theorem 2} Let $D\subset\Bbb C^{\,n}$ be a domain equipped
with $n$ families \thetag{2.6} of nonsingular algebraic curves
algebraically depending on parameters and being in general position.
Then each holomorphic function $f(\bold z)$ in $D$, which is algebraic
in $t_m$ after restriction to any particular curve from any one of
these families, extends to an algebraic function on $\Bbb C^{\,n}$.
\endproclaim
\head
3.Tangent $CR$ fields and the main equations.
\endhead
     First we shall recall briefly some facts of the theory of
$CR$-structures (reader can find more details in \cite{Ch}). Let
$(T^c_pM)_{\Bbb R}$ be the complex tangent space $T^c_pM$ considered
as a vector space over real numbers $\Bbb R$. Then $(T^c_pM)_{\Bbb R}
\,\overset{\Bbb R}\to{\otimes}\,\Bbb C$ is a complexification for
$(T^c_pM)_{\Bbb R}$. Operator $J$ of canonical complex structure
in this complexification is $\Bbb C$\,-linear with respect to
the own complex structure of $T^c_pM$. Therefore it has two
eigenvalues $+i$ and $-i$. Then we have the decomposition
$(T^c_pM)_{\Bbb R}\,\otimes\,\Bbb C=T^c_pM^{1,0}\,\oplus T^c_pM^{0,1}$,
where
$$
T^c_pM^{1,0}=\{(v,-iv): v\in T^c_pM\}\text{ and }
T^c_pM^{0,1}=\{(v,+iv): v\in T^c_pM\}
$$
are the eigenspaces for $J$ corresponding to the eigenvalues $+i$
and $-i$ respectively. The following two maps $v\longmapsto
(v,-iv)$ and $v\longmapsto (v,+iv)$ realize the canonical
isomorphisms $T^c_pM\,\equiv\,T^c_pM^{1,0}$ and $\overline{T^c_pM}
\,\equiv\,T^c_pM^{0,1}$. Thus we have $T^c_pM_{\Bbb R}\,\otimes\,
\Bbb C=T^c_pM\,\oplus\overline{T^c_pM}$. Because of the last
decomposition the sections of the vector-bundles $T^cM$ and
$\overline{T^cM}$ are called the $CR$ vector fields of the types
$(1,0)$ and $(0,1)$ respectively. It is easy to check that the
field $V=\sum_{j=1}^n v_j(z)\partial/\partial z_j$ in $\Bbb C^{\,n}$
is the vector field of the type $(1,0)$ on $M$ if and only if the
vector $(v_1(a),\dots,v_n(a))$ is in $T^c_aM$ for any $a\in M$.\par
     Let us come back to theorem 1. Note that without loss of
generality we can take $p=0$ and $F(p)=0$. Taking $z=(x,y)$ and
$z'=(x',y')$ and denoting $\Bbb C^{\,n}=\Bbb C^{\,k}_x\times
\Bbb C^{\,d}_y$ and $\Bbb C^{\,n'}=\Bbb C^{\,k'}_{x'}\times
\Bbb C^{\,d\,'}_{y'}$ one can bring the defining polynomials for
$M$ and $M'$ to the form
$$
\aligned
p_j&=y_j+\bar y_j + o(|z|),\quad j=1,\dots,d\\
p'_j&=y'_j+\bar y'_j + o(|z|),\quad j=1,\dots,d\,'
\endaligned\tag3.1
$$
For this choice of coordinates we have $T^c_0M=\Bbb C^{\,k}_x=
\{(x,y): y=0\}$ and $T^c_0M'=\Bbb C^{\,k}_{x'}=\{(x',y'):
y'=0\}$ Now let us consider vector fields $T_q$, $q=1,\dots,k$
of the following form
$$
T_q=\Delta(z,\bar z)\frac{\partial}{\partial x_q}-
\sum^d_{j=1} a_{jq}(z,\bar z)\frac{\partial}{\partial y_j}
\tag3.2
$$
where $\Delta$ is the determinant of the following matrix
$$
\Delta_{sj}=\left(\frac{\partial\rho_s}{\partial y_j}
\right)^{j=1,\dots,d}_{s=1,\dots,d}\tag3.3
$$
Everywhere in this paper we shall obey the following rule for
denoting the matrix elements: lower index outside the right
bracket is the row number and upper index is the column number.
The coefficients $a_{jq}$ in \thetag{3.4} we define as follows
$$
a_{jq}=\sum^d_{s=1}\Delta\,b_{js}\,
\frac{\partial\rho_s}{\partial x_q}\tag3.4
$$
where $b_{js}$ is the matrix inverse to the matrix \thetag{3.3}.
According to the elementary facts from linear algebra the
matrix with the elements $\Delta\, b_{js}$ is an conjugate
matrix for \thetag{3.3} i.e. its elements are the algebraic
cofactors for the elements of transpose of the matrix
\thetag{3.3}. Therefore the coefficients of the vector
fields in \thetag{3.2} are polynomials in $x_i$ and $y_i$.\par
     Clear that the restrictions of the fields $T_q$, $q=1,\dots,
k$ form a base of the bundle $T^cM$ over a neighborhood of the
origin in $M$. Moreover, it is obvious that
$$
\Delta(0)=1, a_{jq}(0)=0, j=1,\dots,d, q=1,\dots,k
\tag3.5
$$
Note that because of $dF_0(T^c_0M)\subset T^c_0M'$ , we get by
\thetag{3.1}
$$
\partial F_j/\partial x_q(0)=0, j=k'+1,\dots,n',q=1,\dots,k
\tag3.6
$$\par
     Together with the fields \thetag{3.2} let us consider
the conjugate fields $\overline{T_q}$. Recall that a $C^1$-function
$h$ defined on an open connected subset $U\subset M$ is called
a $CR$-function if for any $z=(x,y)\in U$ one has
$$
\overline{T_q}h=\Delta\,\frac{\partial h}{\partial \bar x_q} -
\sum^d_{j=1}\bar a_{jq}\,\frac{\partial h}{\partial \bar y_j}=0,
q=1,\dots,k
\tag3.7
$$
These are the tangent Cauchy-Riemann equations. This is well known
that a $C^1$-mapping $F:M\longrightarrow M'$ is a $CR$-mapping
(in the sense of the section 2) if and only if each its component
is a $CR$-function on $M$. Indeed, \thetag{3.7} means that
$\bar\partial h(V^q)=0, q=1,\dots,k$,(where $\{V^q\}^k_{q=1}$ is a
base of $T^c_pM$) for any point $p\in M$. Therefore $\bar\partial
h|{T^c_pM}=0$ and \thetag{3.7} means that the restriction
$dh|{T^c_pM}$ is a $\Bbb C$-linear function for any $p\in M$.
\par
     According to the Boggess-Polking theorem \cite{BP} from the
non-degeneracy of the Levi cone of $M$ at the origin $0\in M$
we derive that the mapping $F$ extends holomorphically to a
wedge with the edge $M$. Using the results of \cite{Su} we obtain
that $F$ extends holomorphically to a neighborhood of the origin in
$\Bbb C^{\,n}$. Therefore everywhere below we take $F$ being
holomorphic in a neighborhood $\Omega\owns 0$ in $\Bbb C^{\,n}$
and suppose $U=M\cap\Omega$.\par
     The condition $F(U)\subset M'$ means that $\rho'_j(F,\bar F)=0$
for $z\in U$ and $j=1,\dots,d\,'$. Applying the tangent operators
\thetag{3.2} to both sides of these equalities we obtain
$$
T_q\rho'_j(F,\bar F)=0, q=1,\dots,k,j=1,\dots,d\,'
\text{ for }z\in U\tag3.8
$$\par
     Now we introduce the vector-function $D(F)$ holomorphic in
$\Omega$
$$
D(F)=\left(\frac{\partial F_1}{\partial z_1},\dots,
\frac{\partial F_{n'}}{\partial z_1},\dots,\frac{\partial F_1}
{\partial z_n},\dots,\frac{\partial F_{n'}}{\partial z_n}\right)
$$
Its values are in $\Bbb C^{\,n n'}$.\par
    Left  hand  sides  of  \thetag{3.8}  can  be  considered   as
polynomials
in $z$, $\bar z$, $F$, $\bar F$ and $\partial F_j/\partial z_s$,
$\partial\bar F_j/\partial z_s$. But since $F$ is holomorphic,
we have $\partial\bar F/\partial z_s=0$. Thus, the following
expressions
$$
\Phi_{qj}(z,\bar z, F, \bar F, D(F))=T_q\rho'_j(F,\bar F)
\tag3.9
$$
are polynomials in $z$, $\bar z$, $F$, $\bar F$, $D(F)$.
By $\tilde 0$ we denote the point $(0,0,0,0,D(F)(0))$ in
$\Bbb C^{\,n}\times\Bbb C^{\,n}\times\Bbb C^{\,n'}\times
\Bbb C^{\,n'}\times\Bbb C^{\,nn'}$.
\proclaim{Lemma 3.1} In the set of functions \thetag{3.9} with
$q=1,\dots,k,j=1,\dots,d\,'$ one can choose a subset $\Phi_1=
\Phi_{q(1)j(1)},\dots,\Phi_{k'}=\Phi_{q(k')j(k')}$ such that
the following Jacobian matrix is of rank $k'$
$$
\left(\frac{\partial\Phi_j}{\partial\bar F_s}(\tilde 0)
\right)^{s=1,\dots,k'}_{j=1,\dots,k'}
$$
\endproclaim\par
     \demo{Proof}: Let us  fix $j$, consider the set of functions
$\Phi_{qj}$, $q=1,\dots,k$ and form the following Jacobi matrix
$$
\left(\frac{\partial\Phi_{qj}}{\partial\bar F_s}(\tilde 0)
\right)^{s=1,\dots,k}_{j=1,\dots,k'}
$$
It follows by \thetag{3.1}, \thetag{3.2}, \thetag{3.5}, \thetag{3.6},
\thetag{3.8}, \thetag{3.9} that this matrix can be rewritten in the
form
$$
\left(\sum^{k'}_{r=1}\frac{\partial^2\rho'_j}{\partial z'_r\,
\partial\bar z'_s}(0)\frac{\partial F_r}{\partial z_q}(0)
\right)^{s=1,\dots,k}_{j=1,\dots,k'}\tag3.11
$$\par
     The condition \thetag{1.5} does not depend on the choice of
scalar product in $\Bbb C^{\,n'}$. Therefore one can take the
hermitian scalar product, which defines Levi operators $L^j_0$,
being canonical in the coordinates \thetag{3.1}. Then each
operator $L^j_p$ in the standard basis $e'_r$, $r=1,\dots, k'$ of
$\Bbb C^{\,k'}=T^c_0M'$ has the matrix of the following form
$$
\left(\frac{\partial^2\rho'_j}{\partial z'_j\,\partial\bar z'_s}
(0)\right)^{s=1,\dots,k'}_{j=1,\dots,k'}
$$
Therefore $q$-th the row of the matrix \thetag{3.11} consists of
the coordinates of the vector $L^j_0(dF_0(e_q))$, where
$e_q$, $q=1,\dots,k$ is the standard basis of $\Bbb C^{\,k'}=
T^c_0M'$, i.e.
$$
{}^t(\partial\Phi_{qj}/\partial\bar F_1(\tilde 0),\dots,
\partial\Phi_{qj}/\partial\bar F_k(\tilde 0))=L_0^j(dF_0(e_q))
$$
where ${}^t(\ )$ denotes the transpose matrix. According to
\thetag{2.5} the rank of the set of vectors $L_0^j(dF_0(e_q))$,
$j=1,\dots,d\,'$, $q=1,\dots,k$, is equal to $k$. Thus we
get the proof of the lemma.\qed\enddemo\par
    Let $\Phi_j$, $j=1,\dots,k'$ be the functions chosen according
to lemma 3.1. For $z\in U$ we have
$$
\aligned
&\Phi_j(z,\bar z,F,\bar F,D(F))=0, j=1,\dots,k'\\
&\rho'_s(F,\bar F)=0, s=1,\dots,d'
\endaligned
\tag3.12
$$
\proclaim{Lemma 3.2} The rank of the Jacobian matrix for \thetag{3.12}
with respect to $\bar F$ at the point $\tilde 0$ is equal to $n'$.
\endproclaim
     \demo{Proof} The Jacobi matrix, which is mentioned in lemma, has
the following form
$$
\left(
\vcenter{\hsize 2.0cm
\hbox to 1.9cm{\vsize 1.2cm\hbox to 1.2cm{
\hss$\displaystyle\frac{\partial\Phi_j(\tilde 0)}
             {\partial\bar F_s}$\hss}
\vrule height 0.6cm depth 0.6cm
\hbox to 0.6cm{\hss *\hss}}
\hrule width \hsize
\hbox to 1.9cm{\vsize 1cm\hbox to 1.2cm{\hss$\tilde 0^{k'}_{d'}$\hss}
\vrule height 0.6cm depth 0.4cm
\hbox to 0.6cm{\hss $I_{d'}$\hss}}
}
\right)\tag3.13
$$
Upper left block of the matrix \thetag{3.13} is formed by the matrix
\thetag{3.10}, $\bold 0$ is zero rectangular matrix $d\,'\times k'$
and $I$ is a square unit matrix $d\,'\times k'$. Therefore
because of lemma 3.1 the rank of matrix \thetag{3.13} is equal
to $k'+d\,'=n'$.\qed\enddemo\par
     Applying the complex conjugation to the first $k'$ equations
\thetag{3.13}, we get
$$
\aligned
&P_j(z,\bar z,F,\bar F,\overline{D(F)})=0,\ j=1,\dots,k'\\
&\rho'_s(F,\bar F)=0,\ s=1,\dots,d\,'
\endaligned\tag3.14
$$
This system of the equations is of crucial importance in what follows.
\par
\head
4. Geometry of Segre surfaces.
\endhead
For a real algebraic manifold $M$ of the form \thetag{2.1} the Segre
surface associated to a fixed point $z\in\Bbb C^{\,n}$ is a complex
algebraic set in $\Bbb C^{\,n}$ of the form
$$
Q(z)=\{w\in\Bbb C^{\,n}:\ \rho_j(w,\bar z)=0,\
j=1,\dots,d\}
\tag4.1
$$
We denote by $A$ the graph of the mapping $F$ over a neighborhood
$\Omega \owns 0$ in $\Bbb C^{\,n}$. Let also
$$
A_{\zeta}=\{(z,z')\in\Bbb C^{\,n+n'}:\ z'=F(z),\rho_j(z,\bar\zeta)=0,
\ j = 1,\dots,d\}
\tag4.2
$$
denote the graph of the restriction of $F$ to the Segre surface
$Q(\zeta)$. Evidently, every $A_\zeta$ is a $k$-dimensional complex
manifold ($k=n-d$) in $\Omega\times\Bbb C^{\,n'}$.\par
\proclaim{Lemma 4.1} For any point $\zeta\in\Omega$ the complex
manifold $A_\zeta$ is a piece of a complex $p$-dimensional
algebraic variety $\tilde A_\zeta$ in $\Bbb C^{\,n+n'}$.
\endproclaim
\demo{Proof} It follows by lemma 3.2 that one can apply the
implicit function theorem to the system \thetag{3.14}. We get
$F(z)=R(z,\bar z)$ for $z\in M\cap\Omega$ (where $R$ is a real
analytic function in $\Omega$ algebraic in $z$). By \thetag{3.1}
and the implicit function theorem we get $M\cap\Omega=\{z =(x,y)
\in\Omega:\ y=\phi(x,\bar z)\}$. Therefore,
$$
F(x,\phi(x,\bar z))=R(x,\phi(x,\bar z),\bar z),\tag4.3
$$
for $z=(x,y)\in M\cap\Omega$. Consider antiholomorphic functions
$F^\star(\theta,\xi)=F(\bar\theta,\phi(\bar\theta,\bar\xi))$ and
$R^\star(\theta,\xi)=R(\bar\theta,\phi(\bar\theta,\bar\xi),\bar
\xi)$, where $\theta\in\Bbb C^{\,k}$, $\xi\in\Bbb C^{\,n}$. Then
\thetag{4.3} means that these functions coincide on the manifold
$$
\hat M=\{(\theta,\xi):\ \bar\theta=(\xi_{1},\dots,\xi_{k}),\
\xi\in M\},
$$
that obviously is generic in a neighborhood of the origin in
$\Bbb C^{\,n+k}$. Now it follows by the uniqueness theorem
\cite{P} that $F(\bar\theta,\phi(\bar\theta,\bar\xi))\equiv
R(\bar\theta,\phi(\bar\theta,\bar\xi),\bar\xi)$ in a
neighborhood of the origin in $\Bbb C^{\,n+k}$. Hence,
$$
F(x,\phi(x,\bar\zeta))=R(x,\phi(x,\bar\zeta),\bar\zeta)
$$
for any fixed $\zeta$ from a neighborhood of the origin in
$\Bbb C^{\,n}$. But the following set
$$
\{(x,\phi(x,\bar \zeta)):\ x\text{ runs over a
neighborhood of the origin}\}
$$
coincide with the Segre surface $Q(\zeta)=\{z:\rho_j(z,
\bar\zeta)=0,\ j=1,\dots,d\}$. Thus we get $F(z)=R(z,\bar\zeta)$
for $z\in Q(\zeta)$. Since $R$ was obtained by \thetag{3.14},
the set
$$
A_\zeta=\{(z,z'):\ z\in Q(\zeta),\ z'=F(z)\}
$$
is contained in $(n-d)$-dimensional complex algebraic manifold of
the form
$$
\aligned
&P_j(z,\bar \zeta,z',\overline F(\zeta),\overline {DF(\zeta)})=0,
\ j = 1,\dots,k',\\
&\rho'_s(z',\overline F(\zeta))=0,\ s = 1,\dots,d',\\
&\rho_l(z,\bar \zeta)=0,\ l =1,\dots,d,
\endaligned
\tag4.5
$$
that proves the desired assertion.\qed\enddemo\par
     Fix $\theta\in\Bbb C^{\,k}$ and consider the $d$-parametric
family of the Segre surfaces $Q(\theta,\tau)$, $\tau\in
\Bbb C^{\,d}$.
\proclaim{Lemma 4.2} There exists a neighborhood $U\owns 0$ in
$\Bbb C^{\,n}$ of the form $U=U_{x}\times U_{y}$, $U_{x}\subset
\Bbb C^{\,k}$, $U_{y}\subset\Bbb C^{\,d}$ such that for any fixed
$\theta\in U_x$, the family of Segre surfaces $Q(\theta,\tau)$,
$\tau\in U_y$ has the following properties:
\roster
\item for any $\tau',\tau'' \in U_{y}$ the intersection
      $Q(\theta,\tau')\cap Q(\theta,\tau'')$ is empty;
\item for any $z=(x,y) \in U$ there exists the unique
      $\tau \in U_{y}$ such that $(x,y) \in Q(\theta,\tau)$.
\endroster
\endproclaim
\demo{Proof} One can represent the Segre surface as $Q(\theta,
\tau)=\{z\in U:\ \bar\tau=S(z,\bar\theta)\}$, where $S$ is an
analytic function and $U$ is a neighborhood of the origin. Now,
if $z=(x,y)$ is in $Q(\theta,\tau')\cap Q(\theta,\tau'')$, then
$\bar\tau'=S(z,\bar\theta)=\bar\tau''$; this implies \therosteritem{1}.
For $z=(x,y)$ set $\tau=\bar S(z,\bar\theta)$. Then $z\in Q(\theta,
\tau)$ and we get \therosteritem{2}.
\qed\enddemo\par
     By the implicit function theorem
$$
Q(\theta,\tau)\cap U=\{(x,y)\in U:\ y=R(x,\bar\theta,\bar\tau)\},
$$
where $R$ is an algebraic function, i.e. locally $Q(\theta,\tau)$
is the graph over the coordinate plane $\Bbb C^{\,k}_x=
\Bbb C^{\,k}_{z_1\dots z_k}$. Let $X_j(\bar\theta,\bar\tau)$ be
holomorphic vector fields on $Q(\theta,\tau)$ being the natural
liftings to $Q(\theta,\tau)$ of the coordinate vector fields
$\partial/\partial z_j$, $j=1,\dots,k$ in $\Bbb C^k_{z_1\dots z_k}$.
It follows by lemma 4.2 that for any point $(x,y)\in U$ there exists
the unique surface $Q(\theta,\tau)$, $\tau=S(z,\bar\theta)$ passing
through $(x,y)$. Hence, the fields $Y_j(\bar \theta)=X_j(\bar\theta,
S(z,\bar \theta))$ are holomorphic vector fields in $U$. Their
integral curves evidently are linear sections of the Segre surfaces
by parallel planes and, therefore, form the families of complex
algebraic curves in $\Bbb C^{\,n}$ algebraically depending on
parameters.
\proclaim {Lemma 4.3} The set of the vectors $Y_j(\bar\theta)(0)$,
$j=1,\dots,k$, $\theta$ runs over a neighborhood of the origin in
$\Bbb C^{\,k}$, spans $\Bbb C^{\,n}$.
\endproclaim\par
\demo{Proof} If $(\theta,\tau) \in Q(0)$, i.e. $\rho_j(\theta,\tau,
0,0)=0$, $j=1,\dots,d$, then it follows by \thetag{3.1} and the
implicit function theorem that $\tau=o(|\theta|)$. Now let
$0 \in Q(\theta,\tau)$ (recall that this is equivalent to $(\theta,\tau)
\in Q(0)$). By the implicit function theorem we get
$$
\aligned
&Q(\theta,\tau)=\{(x,y):\ y+\bar\tau=\phi(x,\bar\theta,
y,\bar\tau)\}=\\
&\{(x,y):\ y+\bar\tau=\psi(x,\bar\theta,\bar\tau)\}.
\endaligned
$$
Here
$$
\psi=o(|\theta|)+\bigl<L(x),\ \theta+o(|\theta|)\bigr>+o(|x|),
$$
where
$$
\bigl<L(\xi),\ \eta\bigr>=(\bigl<L_1(\xi),\ \eta\bigr>,\dots,\bigl<
L_d(\xi),\ \eta\bigr>)
$$
is the Levi form of $M$. Hence
$$
Y_j(\theta)(0)=(e_j,\bigl<L(e_j),\ \theta\bigr>+o(|\theta|)) =
(0,\dots,1,\dots,0,\bigl<L(e_j),\ \theta\bigr>+o(|\theta|)),
$$
where $1$ is on the $j$-th position and $e_j$, $j=1,\dots,k$ is the
standard basis of $\Bbb C^{\,k}$.\par
     Assume there exists $\alpha\in\Bbb C^{\,n}\backslash\{0\}$
such that
$$
\bigl<\alpha,Y_j(\bar \theta)(0)\bigr>=0,\ j=1,\dots,k
$$
for any $\theta\in U_x$. Then
$$
\aligned
&\bigl<\alpha,\ Y_j(\bar\theta)(0)\bigr>=\alpha_j+\sum^d_{\nu=1}
(\bigl<L_\nu(e_j),\ \theta\bigr>+o(|\theta|))\alpha_{k+\nu}=\\
&=\alpha_j+\sum^d_{\nu=1}
\alpha_{k+\nu}\bigl<L_{\nu}(e_j),\ \theta\bigr>+o(|\theta|)=\\
&=\alpha_j+\bigl<\sum^d_{\nu=1}\alpha_{k+\nu}L_{\nu}(e_j),\
\theta\bigr>+o(|\theta|)\equiv 0
\endaligned
$$
as a function in $\theta\in\Bbb C^{\,k}$. Therefore $\alpha_j=0$ and
$\sum^d_{\nu=1} \alpha_{k+\nu}L_{\nu}(e_j) = 0$ for $j=1,\dots,k$. This
means that the Levi operators of $M$ are linearly dependent. We get a
contradiction with the condition of the non-degeneracy of the Levi cone
of $M$.\qed\enddemo\par
     Thus, we get $n$ non-singular families of algebraic curves,
algebraically depending on the parameters, in general position near
the origin and the restriction of $F$ on each curve is algebraic.
Now it follows by theorem 2 that $F$ extends to an algebraic mapping
on all $\Bbb C^{\,n}$. This completes the proof of theorem 1 provided
theorem 2 holds.
\head
5. Proof of the theorem 2.
\endhead
     First consider only the $m$-th family of algebraic curves
\thetag{2.6}. Because of nonsingularity of this family we can
treat $t_m, c^{(m)}_1,\dots,c^{(m)}_{n-1}$ in \thetag{2.6} as
new local coordinates in the domain $D$. Let us define the
transformation $\varphi_m(\tau)$ as a translation along the
$t_m$-axis in new curvilinear coordinates
$$
\varphi^{(m)}(\tau): t_m,c^{(m)}_1,\dots,c^{(m)}_{n-1}
\longrightarrow t_m+\tau,c^{(m)}_1,\dots,c^{(m)}_{n-1}
\tag5.1
$$
The transformations $\varphi^{(m)}(\tau)$ form the local one-parameter
group of transformations determined by the vector field of tangent
vectors to the curves of $m$-th family. In the original variables
transformations \thetag{5.1} are given by $n$ algebraic functions
with $n+1$ arguments
$$
\aligned
&\tilde z_1=\varphi_1^{(m)}(\tau,z^1,\dots,z_n)\\
&\text{\hbox to 5cm{\leaders\hbox to 1em{\hss.\hss}\hfill}}\\
&\tilde z_n=\varphi_n{(m)}(\tau,z^1,\dots,z_n)
\endaligned\tag5.2
$$
The algebraicity of the functions $\varphi_i^{(m)}$ in \thetag{5.2}
is a consequence of algebraicity of the curves of $m$-th family
and of their algebraic dependence on parameters in \thetag{2.6}.
\par
     Using the transformations $z\longmapsto\tilde z=\varphi^{(m)}
(\tau) z$ of the form \thetag{5.2} we introduce new local
holomorphic coordinates $t_1,\dots,t_n$ in $D$ as follows
$$
z=\varphi^{(n)}(t_n)\circ\ldots\circ\varphi^{(1)}(t_1) z^0
\tag5.3
$$
where $z^0$ is a fixed point in whose neighborhood these coordinates
are defined (recall that our families of curves are in general position).
The transformation from $z_1,\dots,z_n$ to $t_1,\dots,t_n$ and the
inverse are algebraic.\par
     Note that part of coordinate lines in the local coordinates
$t_1,\dots,t_n$ coincide with the curves of the above families.
Let $f(t_1,\dots,t_n)$ be the representation of the function
$f(z)$ from theorem 2 in the local coordinates $t_1,\dots,t_n$.
Then the following functions
$$
\aligned
&f_1=f(t,0,\dots,0)\\
&f_2=f(t_1,t,0,\dots,0)\\
&\text{\hbox to 5cm{\leaders\hbox to 1em{\hss.\hss}\hfill}}\\
&f_n=f(t_1,\dots,t_{n-1},t)
\endaligned\tag5.4
$$
are algebraic in $t$ and holomorphic in other arguments.
In order to prove theorem 2 it suffices to show the algebraicity
of these functions in all their arguments. We shall proceed by
induction on i (the number of the function $f_i$ in
\thetag{5.4}). However, first we need some preliminaries. \par
     Let $\alpha=(\alpha_1,\dots,\alpha_n)$ be an entire
multiindex and $|\alpha|=\alpha_1+\ldots+\alpha_n$. We denote by
$f_\alpha(z)$ the following derivative
$$
f_\alpha=\frac{\partial^{|\alpha|}f}{\partial z_1^{\alpha_1}\ldots
\partial z_n^{\alpha_n}}\tag5.5
$$
\proclaim{Lemma 5.1} In the assumptions of theorem 2 one can find the
smaller subdomain $D'\subset D$ such that all derivatives $f_\alpha(z)$
are algebraic in $t_m$ after restriction to each curve of any family.
\endproclaim\par
\demo{Proof} The family of curves \thetag{2.6} is non-singular, therefore
the transformation \thetag{2.6} from $z_1,\dots,z_n$ to $t_m,
c^{(m)}_1,\dots,c^{(m)}_{n-1}$ and its back transformation are
implemented by algebraic functions. Hence in place of \thetag{5.5}
we can consider the following derivatives
$$
f_\alpha=\frac{\partial^{\alpha_1+\ldots\alpha_{n-1}}f}
{\partial {c^{(m)}_1}^{\alpha_1}\ldots
 \partial {c^{(m)}_{n-1}}^{\alpha_{n-1}}}\tag5.6
$$
and prove their algebraicity in $t_m$. Differentiation by $t_m$
and transformation to the original variables $z_1,\dots,z_n$ do
not destroy their algebraicity.\par
     In order to prove the algebraicity of the derivatives \thetag{2.6}
we shall use the algebraicity of the function $f(t_m,c^{(m)}_1,\dots,
c^{(m)}_{n-1})$ in $t_m$ for the fixed values of the other
arguments. This means that we have the irreducible polynomial
with unit content in the ring $\Bbb C[f,t]$ such that $f(t_m)$
satisfies the following equation
$$
P(f(t_m),t_m)\equiv 0\tag5.7
$$
(see \cite{VW, L}). Note that the coefficients of the polynomial
\thetag{5.7} and even its degrees in $f$ and $t$ depend on the
parameters $c^{(m)}_1,\dots,c^{(m)}_{n-1})$. Let us define the
following sets being the subsets in the range of values of these
parameters
$$
C_{qk}=\{c^{(m)}_1,\dots,c^{(m)}_{n-1}): \deg_fP=q, \deg_tP=k\}
\tag5.8
$$
The union of the countable number of sets $C_{qk}$ coincides
with the whole range of values of the parameters $c^{(m)}_1,
\dots,c^{(m)}_{n-1}$. This allows us to use the following
well-known Bair theorem. \par
\proclaim{Bair's Theorem} A compete metric space cannot be a countable
union of nowhere dense subsets.
\endproclaim
The proof can be found in  \cite{RS}. We apply this fact to the
of range of parameters  $c^{(m)}_1,\dots,c^{(m)}_{n-1}$, and conclude that
the closure of at least one of the sets  $C_{qk}$ has the non-empty
interior. Choose a domain $D'$ whose natural projection lies
in the interior of such   $C_{qk}$. Also, let us choose in $D'$
a curve of the family \thetag{2.6} with parameters in  $C_{qk}$.
Without loss of generality one can assume that this curve corresponds to
the parameters $c^{(m)}_i=0$ and the point
$t_m=0$ on this curve is in $D'$. For the polynomial
\thetag{5.7} we have
$$
P(f,t)=\sum^q_{i=0}\sum^k_{j=0} a_{ij}\,f^i\,t^j\tag5.9
$$
We normalize the polynomial \thetag{5.9} by setting some of its nonzero
coefficient $a_{rs}$ to be equal to $1$. This polynomial vanishes after
the substitution $f=f(t_m)$ and $t=t_m$. Let us consider the functions
$$
\varphi_{ij}=f(t)^i\,t^j\text{, where }
i=0,\dots,q\text{ and }j=0,\dots,k\tag5.10
$$
They are algebraic in $t$ and depend holomorphically on the
parameters $c^{(m)}_i$. If these parameters vanish, these
functions (as the functions in $t$) are linearly dependent. But
the elimination of the function $\varphi_{rs}$ with $a_{rs}=1$
makes the rest functions linearly independent. Otherwise we would have
another nonzero polynomial $\tilde P(f,t)$ of the form \thetag{5.9}
for which the equality \thetag{5.7} holds. Since $P$ is irreducible,
we have $\tilde P(f,t)=uP(f,t)$ where $u\in\Bbb C[f,t]$. But $\deg_f
\tilde P\le deg_f P$ and $\deg_t\tilde P\le deg_t P$, therefore
$u\in\Bbb C\subset\Bbb C[f,t]$. Comparing the coefficients
$\tilde a_{rs}=0$ and $a_{rs}=1$ we find that the equality
$\tilde P(f,t)=uP(f,t)$ cannot be true for $u\neq 0$.\par
     We denote by $X$ the set of all functions in \thetag{5.10},
and by $X'$ this set without $\varphi_{rs}$. Let us consider the
Taylor expansions in $t$ of the functions \thetag{5.10}. One can
treat their coefficients as infinitely-dimensional vectors (columns)
of the linear space $\Bbb C^\infty$. Such vectors corresponding to
the function from $X'$ form the $\infty\times N$-matrix $A$,
where $N=\#X'$ is the number of functions in the set $X'$. The
columns of $A$ are linearly independent if $c^{(m)}_i=0$. Therefore,
there is a $N\times N$-submatrix $\tilde A$ of $A$ with non-zero
determinant (minor). This minor is holomorphic in $c^{(m)}_i$ and,
therefore, does not vanish in a neighborhood of the origin.
Hence, the columns of $A$ and the functions from $X'$ are linearly
independent for $c^{(m)}_i$ in a neighborhood of the origin.\par
     Let us add the last column $B$ corresponding to the function
$\varphi_{rs}$ to the matrix $A$, and consider the minors of order
$(N+1)$ of the extended matrix  $A|B$. They vanish for $c^{(m)}_i=0$
and for the parameters from the dense set $C_{qk}$. Therefore, they
vanish identically. Thus, the functions from  $X'$ are linearly
independent and the functions of $X$ are linearly dependent for
every $c^{(m)}_i$ in a neighborhood of the origin.
Thus, $\varphi_{rs}$ is a linear combination of the functions from
$X'$. Its coefficients up to the sign coincide with the coefficients of
the polynomial \thetag{5.9}.
They are defined uniquely by linear system with the extended matrix
$(\tilde A|\tilde B)$, where the $N$-column $\tilde B$ is formed
by the elements of $B$ lying on the rows defining $\tilde A$.
Thus, the coefficients of the polynomial \thetag{5.9} are holomorphic
in $c^{(m)}_i$ on a neighborhood of the origin. Now one can
differentiate the equation \thetag{5.7} in $c^{(m)}_i$ and conclude the
proof.\qed\enddemo\par
     Let us consider the functions \thetag{5.4}. One can shrink $D$ to
$D'\subset D$ following to lemma 5.1. Also, one can assume that the
degrees of the polynomials \thetag{5.7} in $f$ and $t_m$ depend only on
$m$ in $D'$. Choose the point $z^0$ from  \thetag{5.3} in a domain
$D'$. This determines the functions \thetag{5.4}. For the function $f_1$
lemma 5.1 gives the algebraicity in $t$ of the derivatives
$$
\frac{\partial^s f_1(t,0,\dots,0)}{\partial t_2^{\,s}}
\tag5.11
$$
The derivatives \thetag{5.11} coincide with the derivatives of the function
$f_2$
from  \thetag{5.4} for $t=0$. In fact,
$$
\left.\frac{\partial^s f_2(t_1,t,0,\dots,0)}{\partial t^s}
\right|_{t=0}=\frac{\partial^s f_1(t_1,0,\dots,0)}
{\partial t_2^{\,s}}\tag5.12
$$
We need the following
\proclaim{Lemma 5.2} An algebraic function $f(t)$ is defined uniquely
by its value and the values of a finite number of its derivatives in a
regular point. If these values depend algebraically on a parameter $\tau$,
then $f=f(t,\tau)$ is an algebraic function function in both variables
$t$ and $\tau$.
\endproclaim
\par Assume the defining irreducible polynomial of the  algebraic
function $f(t)$ has the form \thetag{5.9}. Repeating the above
arguments, we again consider the functions \thetag{5.10} and
their Taylor expansions at a regular point (one can assume it to be
$t=0$). The coefficients of these expansions depends linearly on
$f$ and its derivatives in $t=0$. Considering the non-degenerate
submatrix $\tilde A$, we apply the Cramer rule to the  system
\thetag{5.11} and get the coefficients of the polynomial
\thetag{5.9}. In the second hypothesis of our claim they are
algebraic in $\tau$. By the separate algebraicity principle we
complete the proof of lemma 5.2.\par
     We apply Lemma 5.2 to the function  $f_2(t_1,t,0,\dots,0)$, taking
into account its algebraicity in $t$ and the algebraicity of derivatives
\thetag{5.12} in $t_1$. Therefore, the function $f_2(t_1,t,0,\dots,0)$
is algebraic in both variables. This is the base of the induction.\par
     Assume that the functions  $f_1,\dots, f_m$ in \thetag{5.4} are
algebraic. It follows by lemma 5.1 that the derivatives
$$
\frac{\partial^s f_m(t_1,\dots,t_m,0,\dots,0)}
{\partial t_{m+1}^{\,s}}=\left.
\frac{\partial^s f_{m+1}(t_1,\dots,t_{m+1},0,\dots,0)}
{\partial t_{m+1}^{\,s}}\right|_{t_{m+1}=0}
\tag5.13
$$
are algebraic as well. Lemma 5.2 and the algebraicity of the derivatives
\thetag{5.13} in $t_1,\dots,t_m$ give the induction step from $m$ to
$m+1$. This completes the proof of Theorem 2.

\enddocument
\end